\documentclass{amsart}

\title{The Erd\H os paradox}
\author{Melvyn B. Nathanson}
\address{Department of Mathematics,
Lehman College (CUNY),
Bronx, New York 10468, and 
CUNY Graduate Center, New York, NY 10016}
\email{melvyn.nathanson@lehman.cuny.edu}

\begin{document}

\maketitle

\section*{Prologue}

 The great Hungarian mathematician Paul Erd\H os was born in Budapest 
 on March 26, 1913.  He died alone in a hospital room in Warsaw, Poland, 
 on Friday afternoon, September 20, 1996.  It was sad and ironic that he was alone, 
 because he probably had more friends in more places than any mathematician 
 in the world.  
He was in Warsaw for a conference.   
 Vera S\' os had also been there, but had gone to Budapest on Thursday 
 and intended to return on Saturday with Andr\' as S\' ark\" ozy to travel 
 with Paul to a number theory meeting in Vilnius.  
 On Thursday night Erd\H os felt ill and called the desk in his hotel.  
 He was having a heart attack and was taken to a hospital, 
 where he died about 12 hours later.  No one knew he was in the hospital.  
 When Paul did not appear at the meeting on Friday morning, 
 one of the Polish mathematicians called the hotel.  
 He did not get through, and no one tried to telephone the hotel again for several hours.  
 By the time it was  learned  that Paul was in the hospital, he was dead.  
 
 Vera was informed by telephone on Friday afternoon that Paul had died.  
 She returned to Warsaw on Saturday.  It was decided that Paul should be cremated.  
 This was contrary to Jewish law, but Paul was not an observant Jew 
 and it is not known what he would have wanted.  
 Nor was he buried promptly in accordance with Jewish tradition.  
 Instead, four weeks later, on October 18, there was a secular funeral service 
 in Budapest, and his ashes were buried in the Jewish cemetery in Budapest.  
 
 Erd\H os strongly identified with Hungary and with Judaism.  
 He was not religious, but he visited Israel often, and established a mathematics prize 
 and a post-doctoral fellowship  there.  
 He also established a prize and a lectureship in Hungary.  
 He told me that he was happy whenever someone proved a beautiful theorem, 
 but that he was especially happy if the person who proved the theorem was 
 Hungarian or Jewish.  
 
 Mathematicians from the United States, Israel, and many European countries travelled to Hungary to attend Erdos's funeral.  The following day a conference, 
 entitled ``Paul Erd\H os and his Mathematics,''  took place at the Hungarian Academy of Sciences in Budapest, and mathematicians who were present for the funeral were asked to lecture on different parts of Erd\H os's work.  I was asked to chair one of the sessions, and to begin with some personal remarks about my relationship with Erd\H os and his life and style.  
 
 This paper is in two parts.   
 The first is the verbatim text of my remarks at the Erd\H os memorial conference 
 in Budapest on October 19, 1996.  
 A few months after the funeral and conference I returned to Europe to lecture in Germany.  At Bielefeld someone told me that my eulogy had generated controversy, and indeed, I heard the same report a few weeks later when I was back in the United States. 
 Eighteen years later, on the 100th anniversary of his birth, 
 it is fitting to reconsider  Erd\H os's life and work.

 \section{Eulogy, delivered in Budapest on October 19, 1996}
 
 I knew Erd\H os for 25 years, half my life, but still not very long compared to many people in this room.  His memory was much better than mine; he often reminded me that we proved the theorems in our first paper in 1972 in a car as we drove back to Southern Illinois University in Carbondale after a meeting of the Illinois Number Theory Conference in Normal, Illinois.   He visited me often in Carbondale, and even more often after I moved to New Jersey.  He would frequently leave his winter coat in my house when he left for Europe in the spring, and retrieve it when he returned in the fall.  I still have a carton of his belongings in my attic.  My children Becky and Alex, who are five and seven years old, would ask, ``When is Paul coming to visit again?''  They liked his silly tricks for kids, like dropping a coin and catching it before it hit the floor.  He was tolerant of the dietary rules in my house, which meant, for example, no milk in his espresso if we had just eaten meat.  
 
 He was tough.  ``No illegal thinking,'' he would say when we were working together.  This meant no thinking about mathematical problems other than the ones we were working on at that time.  In other words, he knew how to enforce party discipline.  
 
 Erd\H os loved to discuss politics, especially Sam and Joe, which, in his idiosyncratic language, meant the United States (Uncle Sam) and the Soviet Union (Joseph Stalin).   His politics seemed to me to be the politics of the 30's, much to the left of my own.  He embraced a kind of naive and altruistic socialism that I associate with idealistic intellectuals of his generation.  He never wanted to believe what I told him about the Soviet Union as an ``evil empire.''  I think he was genuinely saddened by the fact that the demise of communism in the Soviet Union meant the failure of certain dreams and principles that were important to him.  
 
 Erd\H os's cultural interests were narrowly focused.  When he was in my house he always wanted to hear ``noise'' (that is, music), especially Bach.  He loved to quote Hungarian poetry (in translation).  I assume that when he was young he read literature (he was amazed that Anatole France is a forgotten literary figure today), but I don't think he read much anymore.
 
 I subscribe to many political journals.  When he came to my  house he would look for the latest issue of \emph{Foreign Affairs}, but usually disagreed with the contents.  Not long ago, an American historian at Pacific Lutheran University published a book entitled \emph{Ordinary Men},\footnote{Christopher R. Browning, \emph{Ordinary Men}, HarperCollins Publishers, New York, 1992.} a study of how large numbers of ``ordinary Germans,'' not just a few SS, actively and willingly participated in the murder of Jews.  He found the book on my desk and read it, but didn't believe or didn't want to believe it could be true, because it conflicted with his belief in the natural goodness of ordinary men.
 
 He had absolutely no interest in the visual arts.  My wife was a curator at the Museum of Modern Art in New York, and we went with her one day to the museum.  It has the finest collection of modern art in the world, but Paul was bored.  After a few minutes, he went out to the scupture garden and started, as usual, to prove and conjecture.  
 
 Paul's mathematics was like his politics.  He learned mathematics in the 1930's in Hungary and England, and England at that time was a kind of mathematical backwater.  For the rest of his life he concentrated on the fields that he had learned as a boy.  Elementary and analytic number theory, at the level of Landau,  graph theory, set theory, probability theory, and classical analysis.  
 In these fields he was an absolute master, a virtuoso.  
 
 At the same time, it is extraordinary to think of the parts of mathematics he never learned.  Much of contemporary number theory, for example.  In retrospect, probably the greatest number theorist of the 1930's was Hecke, but Erd\H os knew nothing about his work and cared less.  Hardy and Littlewood dominated British number theory when Erd\H os lived in England, but I doubt they understood Hecke.  
 
 There is an essay by Irving Segal\footnote{Irving Segal, ``\emph{Noncommutative Geometry} by Alain Connes (book review),'' \emph{Bull. Amer. Math. Soc.} 33 (1996), 459--465} in the current issue of the \emph{Bulletin of the American Mathematical Society}.  He tells the story of the visit of another great Hungarian mathematician, John von Neumann, to Cambridge in the 1930's.  After his lecture, Hardy remarked, ``Obviously a very intelligent young man.  But was that \emph{mathematics}?''
 
 A few months ago, on his last visit to New Jersey, I was telling Erd\H os something about $p$-adic analysis.  Erd\H os was not interested.  ``You know,'' he said about the $p$-adic numbers, ``they don't really exist.''
 
 Paul never learned algebraic number theory.  He was offended -- actually, he was furious -- when Andr\' e Weil wrote that analytic number theory is good mathematics, but analysis, not number theory.\footnote{Weil wrote, ``\ldots there is a subject in mathematics (it's a perfectly good and valid subject and it's perfectly good and valid mathematics) which is called Analytic Number Theory\ldots.  I would classify it under analysis\ldots.'' (\emph{{\OE}uvres Scientifiques Collected Papers}, Springer-Verlag, New York, 1979, Volume III, p. 280).}   Paul's ``tit for tat'' response was that Andr\' e Weil did good mathematics, but it was algebra, not number theory.  I think Paul was a bit shocked that a problem he did consider number theory, Fermat's Last Theorem, was solved using ideas and methods of Weil and other very sophisticated mathematicians.  
 
 It is idle to speculate about how great a mathematician Erd\H os was, as if one could put together a list of the top 10 or top 100 mathematicians of our century.  His interests were broad, his conjectures, problems, and results profound, and his humanity extraordinary.
 
 He was the ``Bob Hope'' of mathematics, a kind of vaudeville performer who told the same jokes and the same stories a thousand times.  When he was scheduled to give yet another talk, no matter how tired he was, as soon as he was introduced to the audience, the adrenaline (or maybe amphetamine) would release into his system and he would bound onto the stage, full of energy, and do his routine for the 1001st time.
 
 If he were here today, he would be sitting in the first row, half asleep, happy to be in the presence of so many colleagues, collaborators, and friends.
 
 Yitgadal v'yitkadash sh'mei raba.
 
 Y'hei zekronoh l'olam.
 
 May his memory be with us forever.\footnote{I ended my eulogy 
 with a sentence in Aramaic and a sentence in Hebrew.  
 The first is the first line of the Kaddish, the Jewish prayer for the dead.  
 Immediately following the second sentence is its English translation.}

 \section{Reconsideration} 
My brief talk at the Erd\H os conference was not intended for publication.   
Someone asked me for a copy, and it subsequently spread via e-mail.  
Many people who heard me in Budapest or who later read my eulogy told me 
that it helped them remember Paul as a human being, but others clearly 
disliked what I said.  I confess I still don't know what disturbed them so deeply. 
It has less to do with Erd\H os, I think, than with the status of ``Hungarian mathematics'' 
in the scientific world.\footnote{cf. L. Babai, ``In and out of Hungary: 
Paul Erd\H os, his friends, and times,'' 
in: \emph{Combinatorics, Paul Erd\H os is Eighty (Volume 2), Keszthely (Hungary) 1993}, 
Bolyai Society Mathematical Studies, Budapest, 1996, pp. 7--95.}  

Everyone understands that Erd\H os was an extraordinary human being 
and a great mathematician who made major contributions to many parts 
of mathematics.  He was a central figure in the creation of new fields, 
such as probabilistic number theory and random graphs.  
This part of the story is trivial.  

It is also true, understood by almost everyone, and not controversial, 
that Erd\H os did not work in and never learned the central core 
of twentieth century mathematics.  
It is amazing to me how great were Erdos's contributions to mathematics, 
and how little he knew.  
He never learned, for example, the great discoveries in number theory 
that were made at the beginning of the twentieth century.  
These include, for example, Weil's work on diophantine equations, 
Artin's class field theory, and Hecke's monumental contributions 
to modular forms and analytic number theory.  
Erd\H os apparently knew nothing about Lie groups, Riemannian manifolds, 
algebraic geometry, algebraic topology, global analysis, or the deep ocean 
of mathematics connected with quantum mechanics and relativity theory.  
These subjects, already intensely investigated in the 1930's, 
were at the heart of twentieth century mathematics.  
How could a great mathematician not want to study these things?\footnote{This suggests the fundamental question:  How much, or how little, must one know in order to do great mathematics?}
This is the first Erd\H os paradox.  

In the case of the Indian mathematician Ramanujan, whose knowledge was also deep but narrow, there is a discussion in the literature abut the possible sources of his mathematical education.  The explanation of Hardy\footnote{``It was a book of a very different kind, Carr's \emph{Synopsis}, which first aroused Ramanujan's full powers,'' according to G. H. Hardy, in his book \emph{Ramanujan}, Chelsea Publishing, New York, 1959, p. 2} and others is that the  only serious book that was accessible to Ramanujan in India was Carr's \emph{A Synopsis of Elementary Results in Pure and Applied Mathematics}, and that Ramanujan lacked a broad mathematical culture because he did not have access to books and journals in India.  But Hungary was not India; there were libraries, books, and journals in Budapest, and in other places where Erd\H os lived in the 1930's and 1940's.  

For the past half century,  ``Hungarian mathematics'' has been a term of art   
to describe the kind of mathematics that Erd\H os 
did.\footnote{For example, Joel Spencer, ``I felt \ldots I was working on `Hungarian  
mathematics','' quoted in Babai, \emph{op. cit.}}
It includes combinatorics, graph theory, combinatorial set theory, and elementary 
and combinatorial number theory.  Not all Hungarians do this kind of mathematics, 
of course, and many non-Hungarians do Hungarian mathematics.  
It happens that combinatorial reasoning is central to theoretical computer science, 
and ``Hungarian mathematics''  commands vast respect in the computer science world.  
It is also true, however, that for many years combinatorics did not have the highest 
reputation among mathematicians in the ruling subset of the research community, 
exactly because combinatorics was concerned largely with questions 
that they believed (incorrectly) were not central to twentieth century  
mathematics.\footnote{For example, 
S. Mac Lane criticized ``emphasizing too much of a Hungarian view of mathematics,'' 
in: ``The health of mathematics,'' \emph{Math.Intelligencer} 5 (1983), 53--55.
}

In a volume in honor of Erd\H os's 70th birthday, Ernst Straus wrote, 
``In our century, in which mathematics is so strongly dominated 
by `theory constructors' [Erd\H os] has remained the prince of problem solvers 
and the absolute monarch of problem posers.''\footnote{E. G. Straus, 
``Paul Erd\H os at 70,'' \emph{Combinatorica} 3 (1983), 245--246.  
Tim Gowers revisited this notion in his essay, ``The two cultures of mathematics,'' 
published in \emph{Mathematics: Frontiers and Perspectives}, 
American Mathematical Society, 2000.} I disagree.  There is, as Gel'fand often said, 
only one mathematics.  There is no separation of mathematics into ``theory'' 
and ``problems.''  But there is an interesting lurking issue.

In his lifetime, did Erd\H os get the  recognition he deserved? 
Even though Erd\H os received almost every honor that can be given 
to a mathematician, some of his friends believe that he was still insufficiently 
appreciated, and they are bitter on his behalf.  
He was awarded a Wolf Prize and a Cole Prize, 
but he did not get a Fields Medal or a permanent professorship 
at the Institute for Advanced Study.  
He traveled from one university to another across the United States, 
and was never without an invitation to lecture somewhere, but his mathematics 
was not highly regarded by the power brokers of mathematics.   
To them, his methods were insufficiently  abstruse and obscure; 
they did not require complicated machinery.   
Paul invented diabolically clever arguments from arithmetic, combinatorics, 
and probability  to solve problems.  But the technique was too simple,  
too elementary.  It was suspicious.  The work could not be ``deep.''
 
 None of this seemed to matter to Erd\H os, who was content to prove 
 and conjecture and  publish more than 1,500 papers.  
  
 Not because of politicking, but because of computer science and 
 because his mathematics was always beautiful, in the past decade 
 the reputation of Erd\H os and the respect paid to discrete mathematics 
 have increased exponentially.  
 The \emph{Annals of Mathematics} will now publish papers in combinatorics, 
 and the most active seminar at the Institute for Advanced Study is 
 in discrete mathematics and theoretical computer science.  
 Fields Medals are awarded to mathematicians who solve 
 Erd\H os-type problems.  Science  has changed.

 In 1988 Alexander Grothendieck was awarded the Crafoord Prize 
 of the Swedish Academy of Sciences.  
 In the letter to the Swedish Academy in which he declined the prize, he wrote,
``Je suis persuad\' e que la seule \' epreuve d\' ecisive pour la f\' ecundit\' e 
d'id\' ees ou d'une vision nouvelles est celle du temps.  
La f\' econdit\' e se reconnait par la prog\' eniture, 
et non par les honneurs.''\footnote{``I believe that time gives the only definite proof 
of the fertility of new ideas or a new vision.  
We recognize fertility by its offspring, and not by honors.''}
 
 Time has proved the fertility and richness of Erd\H os's work.  
 The second  Erd\H os paradox is that his methods and results, 
 considered marginal in the twentieth century, 
 have become  central in  twenty-first century mathematics.  
 
 May his memory be with us forever.

 \end{document}